\documentclass{article}
\usepackage{amsfonts,amssymb,amscd, amsmath}
\usepackage{latexsym, graphicx}
\usepackage[affil-it]{authblk}
\usepackage{algorithmic}
\usepackage{algorithm}
\usepackage{tikz}

\usepackage{enumerate}

\newtheorem{theorem} {Theorem}[section]
\newtheorem{lemma}[theorem]{Lemma}

\newtheorem{problem}[theorem]{Problem}

\usepackage{xcolor}
\newcommand*\circled[1]{\tikz[baseline=(char.base)]{
		\node[shape=circle,draw,inner sep=2pt] (char) {#1};}}

\begin{document}

\title{On the structure of (dart, odd hole)-free graphs}
\author{Ch\'inh T. Ho\`ang }
\date{}
\maketitle
\begin{center}

Department of Computer Science and Physics\\ 
Wilfrid Laurier University\\ 
Waterloo, Ontario,  Canada, N2L 3C5\\
\texttt{choang@wlu.ca}

\end{center}

\begin{abstract}
A hole is a chordless  cycle with at least four vertices. A hole
is odd if it has an odd number of vertices. A dart is a graph
which vertices $a,b,c,d,e$ and edges $ab,bc,bd,be,cd,de$. 
Dart-free graphs have been actively studied in the literature.
We prove that a (dart, odd hole)-free graph is perfect, 
or does not contain a stable set on
three vertices, or is the join or co-join of two smaller graphs. Using this
structure result, we design a polynomial-time algorithm for
finding an optimal colouring of (dart, odd hole)-free graphs. A graph $G$ is perfectly
divisible if every induced subgraph $H$ of $G$ contains a set $X$
of vertices such that $X$ meets all largest cliques of $H$, and
$X$ induces a perfect graph. The chromatic number of a perfectly
divisible graph $G$ is bounded by $\omega^2$ where $\omega$
denotes the number of vertices in a largest clique of $G$. We
prove that (dart, odd hole)-free graphs are perfectly divisible.
\end{abstract}

\section{Introduction}\label{sec:introduction}
A {\em hole} is a chordless cycle with at least four vertices.  An
{\it antihole} is the complement of a hole. A hole is {\it odd} if
it has an odd number of vertices. A graph is odd-hole-free if it
does not contain, as an induced subgraph, an odd hole.
Odd-hole-free graphs have been studied in connection with perfect graphs
(definitions not given here will be given later.) 
Chudnovsky, Scott, Seymour, and Spirkl (\cite{ChuSey2020}) designed a polynomial-time algorithm for
finding an odd hole if one exists, in a graph. 
A result of Kr\'al , Kratochv\'il,   Tuza, and Woeginger (\cite{KraKra2001}) shows that determining the chromatic 
number of an odd-hole-free graph in NP-complete.

A {\it dart} is the graph with vertices $a,b,c,d,e$ and edges $ab,bc,bd,be,cd,de$. 
Dart-free graphs generalize the well studied class of
claw-free graphs. In this paper, we study (dart, odd hole)-free
graphs. We prove that a (dart, odd hole)-free graph is perfect, 
or does not contain a stable set on
three vertices, or is the join or co-join of two smaller graphs. 
We will use this structural
result to design a polynomial-time algorithm for finding an optimal coloring of
a (dart, odd hole)-free graph. Interesting results have been obtained for dart-free graphs. 
Karthick, Maffray and Pastor (\cite{KarMaf2019}) designed a polynomial-time algorithm to 
(optimally) color a ($P_5$, dart)-free graph. Branst\"adt and Mosca (\cite{BraMos2015})
designed a polynomial-time algorithm to find a largest stable set of a (dart, odd hole)-free graph.
Critical dart-free graphs were studied by Xia, Jooken, Goedgebeur and Huang (\cite{XiaJoo2025}). 
Sales and Reed (\cite{SalMaf2004}) studied ``perfectly contractile'' graphs that are dart-free. Chv\'atal, Fonlupt, Sun and Zemirline (\cite{ChvFon2002}) gave a polynomial-time algorithm to recognize dart-free perfect graphs. 

A graph $G$ with at least one edge  is {\it $k$-divisible}  if the
vertex-set of each of its induced subgraphs $H$ with at least one
edge can be partitioned into $k$ sets, none of which contains a
largest clique of $H$. It is easy to see that the chromatic number
of a $k$-divisible graph is at most $k^{\omega -1}$. It was conjectured
by Ho\`ang and McDiarmid ({\cite{HoaMcd1999}, \cite{HoaMcd2002}})
that every odd-hole-free graph is 2-divisible. They proved the
conjecture for claw-free graphs (\cite{HoaMcd2002}). Ho\`ang (\cite{Hoa2018}) proved the conjecture for banner-free graphs.  
Dong, Song and Xu (\cite{DonXu2022}) proved
the conjecture for dart-free graphs.

A graph $G$ is {\it perfectly divisible} (\cite{Hoa2018})  if every induced subgraph
$H$ of $G$ contains a set $X$ of vertices such that $X$ meets all
largest cliques of $H$, and $X$ induces a perfect graph. The
chromatic number of a perfectly divisible graph $G$ is bounded by
$\omega^2$ where $\omega$ denotes the number of vertices in a
largest clique of $G$. We will prove that (dart, odd hole)-free
graphs are perfectly divisible.

In Section~\ref{sec:background}, we give the definitions used in
this paper and discuss background results that we need to establish our results. 
In Section~\ref{sec:structure}, we prove that (dart, odd
hole)-free graphs are perfectly divisible. We will also give a proof, different from that of \cite{DonXu2022}, 
that (dart, odd hole)-free graphs are 2-divisible.  
In Section~\ref{sec:optimization}, we give polynomial-time algorithms
for finding a minimum coloring of a
(dart, odd hole)-free graph.
Finally, in Section \ref{sec:conclusions}, we discuss open problems related to our work. 
\section{Definitions and background}\label{sec:background}
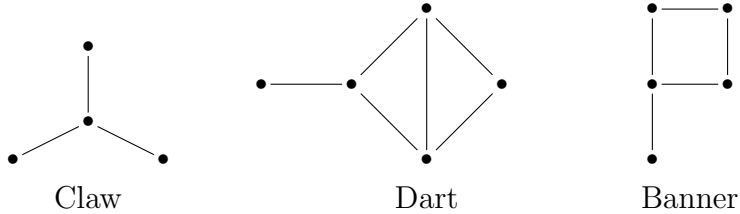
\begin{figure}
\newcommand{\ang}{17}
\newcommand{\sep}{15}
\begin{tikzpicture}[scale=1]
	\tikzstyle{every node}=[font=\small]

\newcommand{\size}{1}

\newcommand{\claw}{1}{
	\path (0,0) coordinate (g1);
	\path (g1) +(\size / 2, \size / 2) node (g1_1){}; 
	\path (g1) +(\size / 2, 1.5 * \size) node	(g1_2){}; 
	\path (g1) +(-\size / 2, 0) node (g1_3){};
	 \path (g1)	+(1.5 * \size, 0) node (g1_4){};
	 
	\foreach \Point in {(g1_1),(g1_2),(g1_3),(g1_4)}{
		\node at \Point {\textbullet};
	}
	\draw  (g1_1) -- (g1_2)
	(g1_1) -- (g1_3)
	(g1_1) -- (g1_4);
	\path (g1) ++(\size  / 2,-\size / 2) node[draw=none,fill=none] { {\large Claw}};
}

\newcommand{\g3}{3}{
	\path( 4.5 * \size, 0) coordinate(g3);
	\path(g3) +(\size / 2, 0) node(g3_1){}; 
	\path(g3) +(-0.5, \size) node(g3_2){};
	\path(g3) +( 1.5 * \size, \size) node(g3_3){};
	\path(g3) +(\size / 2, 2 * \size) node(g3_4){};
	\path(g3) +(-1.7 * \size, \size) node(g3_5){};
	
	\foreach \Point in {(g3_1),(g3_2),(g3_3),(g3_4), (g3_5)}{
		\node at \Point {\textbullet};
	}
	\draw   (g3_1) -- (g3_2)
	(g3_1) -- (g3_3)
	(g3_1) -- (g3_4)
	(g3_4) -- (g3_2)
	(g3_4) -- (g3_3)
	(g3_2) -- (g3_5);
	\path (g3) ++(\size  / 2,-\size / 2) node[draw=none,fill=none] { {\large Dart}};
	
}

\newcommand{\g4}{4} {
	\path( 8 * \size, 0) coordinate(g4);
	\path(g4) + (0,0)  node(g4_1){}; 
	\path(g4) + ( 0 , \size) node(g4_2){};  
	\path(g4) + ( \size , \size) node(g4_3){};  
	\path(g4) + ( 0 , 2 * \size) node(g4_4){};  
	\path(g4) + ( \size , 2*\size) node(g4_5){};  
	
	\foreach \Point in {(g4_1),(g4_2), (g4_3), (g4_4), (g4_5)}{
		\node at \Point {\textbullet};
	}
	\draw   (g4_1) -- (g4_2)
	(g4_2) -- (g4_3) 
	(g4_2) -- (g4_4) 
	(g4_3) -- (g4_5) 
	(g4_4) -- (g4_5);
	\path (g4) ++(\size  / 2,-\size / 2) node[draw=none,fill=none] { {\large Banner}};
}

 \end{tikzpicture}
\caption{The claw, dart and banner}\label{fig:claw-dart}

\end{figure}

The {\it claw}, {\it dart} and {\it banner} are represented in Figure~\ref{fig:claw-dart}.
Let $C_k$ denote the cordless cycle on $k$ vertices. Let $P_k$ denote the chordless path on $k$ vertices.
A {\it hole} is the graph $C_k$ with $k \geq
4$. A hole is {\it odd} if it has an odd number of vertices.  An
{\it antihole} is the complement of a hole. A clique on three
vertices is called a {\it triangle}. For a given graph $H$, it is
customary to let {\it co}-$H$ denote the complement of $H$. Thus,
a co-triangle is the complement of the triangle. Let $L$ be a collection 
of graphs. A graph $G$ is $L$-free if $G$ does not contain
an induced subgraph isomorphic to a graph in $L$. In particular,
a graph is {\it (dart, odd hole)-free} if it does not contain an
induced subgraph isomorphic to a dart or an odd hole.

Let $G$ be a graph. Then $\chi(G)$ denotes the chromatic number of
$G$, and $\omega(G)$ denotes, the {\it clique number} of $G$, that
is, the number of vertices in a largest clique of $G$. The number of vertices in a largest stable set of $G$ is denoted by  
$\alpha(G)$. The {\it girth} of $G$ is the length of the shortest cycle of $G$ (if a cycle exists).

Let $v$ be a vertex of $G$ 
and $X$ be a subset of $V(G)$ that does not contain $v$. 
We say that $v$ is {\it X-complete} if $v$ is adjacent to all vertices 
of $X$, and $v$ is {\it X-anticomplete} if $v$ is non-adjacent to all vertices of $X$.
Let $A$ and $B$ be two disjoint subsets of $V(G)$. We say that
$A$ and $B$ form a {\it join} (respectively, {\it co-join})  if there are all (respectively, no) edges between $A$ and $B$, that is, 
every vertex $x \in A$ is $B$-complete (respectively, $B$-anticomplete). The join of $A$ and $B$ is denoted by $A \;\circled{1}\; B$, and the co-join of $A$ and $B$ is denoted by $A \;\circled{0}\; B$.
%
For a
graph $G$ and a set $X \subseteq V(G)$, $G[X]$ denotes the
subgraph of $G$ induced by $X$.

A graph $G$ is {\it perfect} if $\chi(H) = \omega(H)$ for every
induced subgraph $H$ of $G$. A graph is {\it Berge} if it does not
contain as induced subgraph an odd hole or odd antihole.  The Perfect Graph Theorem, proved by Lov\'asz
(\cite{Lov1972}), states that a graph is perfect if and only if
its complement is. The Strong Perfect Graph Theorem, proved by
Chudnovsky,  Robertson, Seymour, and  Thomas (\cite{ChuRob2006}),
states that a graph is perfect if and only if it is Berge. Both of the above results were long
standing open problems proposed by Berge (\cite{Ber1961}).
Chudnovsky et al. (\cite{ChuCor2005}) designed  an $O(n^9)$
algorithm for recognizing Berge graphs. Gr\"otschel, Lov\'asz and Schrijver
(\cite{GroLov1984}) designed a polynomial-time 
algorithm for finding a largest clique and a minimum coloring of a
perfect graph. Corollary 1 in Kr\'al, Kratochiv\'il, and Tuza (\cite{KraKra2001})
shows that it is NP-hard to compute the chromatic number of an odd-hole-free graph.
For more infomation on perfect graphs, the reader is referred to 
the books of Golumbic (\cite{Gol1980}), Berge and  Chv\'atal (\cite{BerChv1984}), 
Ramirez Alfonsin and Reed (\cite{RamRee2001}), and the survey paper of Ho\`ang and Sritharan (\cite{HoaSri2015}).

\section{The structure of (dart, odd hole)-free graphs}\label{sec:structure}
In this section, we prove a theorem on the structure of (dart, odd hole)-free graphs that we will later 
use to design a polynomial-time algorithm to optimally color such graphs.
\begin{theorem}\label{thm:main}
Let $G$ be a (dart, odd hole)-free graph. Then one of the
following holds for $G$.
\begin{description}
 \item[(i)] $G$ is perfect.
 \item[(ii)] $\alpha(G) \leq 2$.
 \item[(iii)] There is a partition of $V(G)$ into two sets $A$, $B$ such that 
     \begin{description}
     	\item (a) $A \; \circled{0} \; B$, or
     	\item (b) $A \; \circled{1} \; B$, $\alpha(G[A]) = 2$, and $G[B]$ is a perfect graph
     \end{description}

\end{description}
\end{theorem}

To prove Theorem~\ref{thm:main}, we will need the following results that were established by Dong, Song and Xu in \cite{DonXu2022}. In the following two lemmas,  $G$ is a connected (dart, odd hole)-free graph and $A$ is an odd antihole of $G$. 
\begin{lemma}[Lemma 2.2 in \cite{DonXu2022}]\label{lem:co-triangle}
	No vertex of $A$ belongs to a co-triangle of $G$.
\end{lemma}
\begin{lemma}[Lemma 2.4 in \cite{DonXu2022}]\label{lem:join} 
	Every vertex in a co-triangle of $G$ is $A$-complete.
\end{lemma}
We note the proofs of the above two lemmas are simlar to the arguments used by Chv\'atal et al (\cite{ChvFon2002}) to design a recognition algorithm for dart-free perfect graphs.

\noindent {\it Proof of Theorem \ref{thm:main}}. Let $G$ be a (dart, odd hole)-free graph. We may assume that $G$ is connected, 
for otherwise $G$ is the co-join of two graphs and (iii.a) holds. We may assume $G$ contains an odd antihole, for otherwise (i) holds. Now we may assume that $G$ contains a co-triangle, for otherwise (ii) holds. Let $T$ be the set of vertices of $G$ that belongs to a co-triangle. By Lemmas \ref{lem:co-triangle} and \ref{lem:join}, we have $A \cap T = \emptyset$ and $A \; \circled{1} \; T$.
Define $R = V(G) \setminus (A \cup T)$. By definition, vertices in $R$ do not belong to a co-triangle of $G$. Therefore, we have 
\begin{equation}\begin{minipage}{0.8\linewidth}\label{eq:alpha=2}
		$\alpha(G[A \cup R]) = 2$.
\end{minipage}\end{equation}

We claim that 
\begin{equation}\begin{minipage}{0.8\linewidth}\label{eq:total-join}
$T \; \circled{1} \; ( A \cup R)$.
\end{minipage}\end{equation}
Suppose that (\ref{eq:total-join}) failed. So there is a vertex $r \in R$ that is not $T$-complete. Let $t_1$ be a vertex in $T$ that is not adjacent to $r$. By definition, there are vertices $t_2, t_3 \in T$ that together with $t_1$ form a co-triangle of $G$. Since $r$ does not belong to any co-triangle of $G$, $r$ must be adjacent to both $t_2$ and $t_3$. Vertex $r$ must be adjacent to some vertex $a \in A$, for otherwise, $r$ form a co-triangle with some two vertices of $A$, a contradiction. But now, $t_1, a, r_1, t_2, t_3$ induce a dart. So, 
(\ref{eq:total-join}) holds. 

Consider the graph $G[T]$. Each vertex of $G[T]$ belongs to a co-triangle. So, by Lemma \ref{lem:co-triangle}, no odd antihole with at least five vertices contains a vertex of $T$. Thus, $G[T]$ contains no odd antihole, and therefore is perfect by the hypothesis of the theorem. So, (iii.b) holds. $\Box$  

Next, we will prove
\begin{theorem}\label{thm:perfect}
	(Dart, odd hole)-free graphs are perfectly divisible. 
\end{theorem}

We will need the following result. 

\begin{lemma}[Ho\`ang, Lemma 7.3 in \cite{Hoa2018}]\label{lem:alpha-2}
	 Graphs $G$ with $\alpha(G) = 2$ are perfectly divisible 
\end{lemma}	 
\noindent {\it Proof of Theorem \ref{thm:perfect}. }
Let $G$ be a (dart, odd hole)-free graph. We will prove by induction on the number of vertices.  By the induction hypothesis, we may assume every proper induced subgraph of $G$ is perfectly divisible, and we only need to prove that there is a partition of $V(G)$ into two sets $A$, $B$, such that $G[B]$ is perfect and every largest clique of $G$ must contain a vertex of $B$. 

We know that the conclusions of Theorem \ref{thm:main} holds for $G$. If $G$ is perfect, then we are obviously done. If $\alpha(G) \leq 2$, then we are done by Lemma \ref{lem:alpha-2}. Thus there is a partition of $V(G)$ into  two sets of vertices $A$, $B$ that satisfies (iii.a) or (iii.b).  At this point we would be done by Lemma 7.2 in \cite{Hoa2018}, but we can provide a simple proof here. In the case of (iii.b), every largest clique of $G$ contains a vertex of $B$. So $G[B]$ is perfect and meets all largest clique of $G$, and the theorem holds. Consider case (iii.a). By the induction hypothesis, $G[A]$  and $G[B]$ are perfectly divisible. So, $G[A]$ (respectively, $G[B]$) can be partitioned into two sets $A_1$ and $A_2$ (respectively, $B_1$ and $B_2$) such that $G[A_1]$ (respectively, $G[B_1])$)  is perfect and $\omega(G[A_2]) < \omega(G[A])$ (respectively, $\omega(G[B_2]) < \omega(G[B])$). Now we only need to show that $G$ contains a graph $G'$ as induced subgraph  such that $G'$ is perfect and $\omega(G \setminus G') < \omega(G)$. Let $G' = G[A_1 \cup B_1]$. Since both $G[A_1]$ and $G[B_1]$ are perfect and there are no edges between $G[A_1]$ and $G[B_1]$, we know that $G'$ is perfect. Let $K$ be a largest clique of $G$. Since $A \; \circled{0} \; B$, we may assume without loss of generality that $K \cap A \not = \emptyset$. It follows that $K$ lies entirely in $A$ and is met by $A_1$, that is, $\omega(G \setminus G') < \omega(G)$. 
$\Box$

Now, we will use Theorem \ref{thm:main} to give a new proof of the result of \cite{DonXu2022}.
\begin{theorem}[Dong, Song and Xu \cite{DonXu2022}]\label{thm:2-div}
	(Dart, odd hole)-free graphs are 2-divisible
\end{theorem}
We will need a result of Ho\`ang and McDiarmid (Theorem 1 in \cite{HoaMcd2002}). 
\begin{theorem}[Ho\`ang and McDiarmid \cite{HoaMcd2002}]\label{thm:claw}
	(Claw, odd hole)-free graphs are 2-divisible.
\end{theorem}
Let $G$ be a graph and let $H$ be a set of vertices of $G$. The set $H$ is {\it homogeneous} if $2 \leq |H| < |V(G)|$ and every vertex of $V(G) \setminus H$ is $H$-complete or $H$-anticomplete. A graph $G$ is {\it minimally non-2-divisible} if $G$ is not 2-divisible but every proper induced subgraph of $G$ is. We will need the result below (Lemma 3 in \cite{HoaMcd2002}).
\begin{lemma}[Ho\`ang and McDiarmid \cite{HoaMcd2002}]\label{lem:homogeneous}
	No minimally non-2-divisible graph can contain a homogeneous set. 
\end{lemma}  

\noindent {\it Proof of Theorem \ref{thm:2-div}}.   Let $G$ be a (dart, odd hole)-free graph. We may assume $G$ has at least three vertices. We will prove by contradiction. Suppose that $G$ is not 2-divisible. Then $G$ contains an induced subgraph $G'$ that is minimally non-2-divisible. We know the conclusions of Theorem \ref{thm:main} holds for $G'$. If $G'$ is perfect, then $G'$ is 2-divisible, a contradiction. If $\alpha(G) \leq 2$, then $G'$ is claw-free and we have a contradiction by Theorem \ref{thm:claw}. Thus there is a partition of $V(G')$ into  two sets of vertices $A$, $B$ that has $A \; \circled{0} \;B$ or $A \; \circled{1} \; B$. Thus, $A$ or $B$ is a homogeneous set of $G'$, but this is a contradiction to Lemma \ref{lem:homogeneous}. $\Box$

\section{Optimizing (dart, odd hole)-free graphs}\label{sec:optimization}
In this section, we consider the following four optimization
problems: finding a minimum coloring, finding a minimum clique
cover (the minimum coloring in the complement of the graph),
finding a largest stable set, and finding a largest clique. 
These four problems have been much studied in the context of perfect graphs
and related graph classes. 
It is
interesting to note that for our class of (dart, odd hole)-free
graphs, the coloring and largest stable set problems are solvable
in polynomial time, whereas the clique cover and largest clique
problems are NP-hard.
\subsection{Finding a minimum coloring in polynomial time}\label{sec:coloring}
In this section, we describe an algorithm to optimally  color a
(dart, odd hole)-free graph. Using Theorem \ref{thm:main}, we consider the three cases (i), (ii) and (iii). 
If $G$ is perfect, we can use the algorithm of Gr\"otchel et al (\cite{GroLov1984}) to color our graph. 

Concerning graphs $G$ with $\alpha(G) \leq 2$, it is well known
that a minimum coloring of  $G$ can be found by finding a maximum
matching $M$ in the complement $\overline{G}$ of $G$. (Let $t$ be
the number of edges of $M$. Then we have $\chi(G) = n - t$.)  The
algorithm of Micali and Vazirani (\cite{MicVaz1980}) finds a
maximum matching of a graph in $O(m \sqrt{n})$ time.

Now consider the case (iii). There is a partition of $V(G)$ into two sets $A$, $B$ such that $A \; \circled{1} \; B$ or $A \; \circled{0} \; B$. We can recursively compute $\chi(G[A])$ and $\chi(G[B])$.  
If $A \; \circled{0} \; B$, then we have $\chi(G) = max(\chi(G[A]), \chi(G[B]))$. 
If $A \; \circled{1} \; B$, then we have $\chi(G) = \chi(G[A]) + \chi(G[B])$, and we are done.

\subsection{Finding a largest clique is hard}\label{sec:find-clique}
It was proved by Poljak (\cite{Pol1974}) that it is NP-hard to
determine $\alpha(G)$ for a triangle-free graph $G$. We are going to use Poljak's argument to show that 
finding a largest clique of a (dart, odd hole)-free graph is NP-hard. Our argument 
follows that of  \cite{Hoa2018} (Section 5.2 of \cite{Hoa2018}).

Let $G$ be a graph.
Construct a graph $f(G)$ from $G$ by, for each edge $ab$ of $G$,
replacing $ab$ by an induced path on three edges, that is, we
subdivide the edge $ab$ twice. Then it is easy to see that
$\alpha(f(G)) = \alpha(G) + m$, where $m$ is the number of edges
of $G$. The graph $f(G)$ is triangle-free. By repeatedly applying
this construction, we can see that it is NP-hard to compute
$\alpha(G)$ for a graph $G$ of any given minimum girth $g$. In particular,
it is NP-hard to compute $\alpha(G)$ for a (triangle, $C_5$)-free
graph $G$, or equivalently, to compute $\omega(F)$ for a
(co-triangle, $C_5$)-free graph $F$. Since (co-triangle,
$C_5$)-free graphs are (dart, odd hole)-free, it is NP-hard to
compute $\omega(G)$ for a (dart, odd hole)-free graph.
\subsection{Finding a minimum clique cover is hard}\label{sec:find-clique-cover}
Jensen and  Toft (Section 10.3 of
\cite{JenTof1995}) showed that it is NP-complete to decide, for any fixed
integer $g$, whether a graph $G$ of girth at least $g$ is
6-colorable (this result uses the Haj\'os construction for graphs
of high chromatic number). Thus, it is NP-hard to find a minimum
coloring of a $(C_3, C_5)$-free graph, or equivalently, to
find a minimum clique cover of a (co-triangle,  $C_5$)-free
graph $G$. Since such a graph $G$ is (dart, odd hole)-free, it
is NP-hard to find a minimum clique cover of a (dart, odd
hole)-free graph. The argument used here is slightly modified from that of \cite{Hoa2018} (Section 5.3 of \cite{Hoa2018}).
\subsection{Finding a largest stable set in polynomial time}\label{sec:find-stable-set}
Let $\Pi$ denote a (hereditary) graph property. A graph $G$ is {\it nearly $\Pi$} if for all $v \in V(G)$, the subgraph of $G$
induced by the non-neighbors of $G$ (different from $v$) has property $\Pi$, that is, the graph $G[V(G) \setminus (\{v\} \cup 
N(v))]$ has property $\Pi$, where $N(v)$ denotes the set of vertices adjacent to $v$. For example, a graph $G$ is {\it nearly perfect} if for every vertex $v$, the graph $G[V(G) \setminus (\{v\} \cup 
N(v))]$ is perfect. The notion of nearly $\Pi$ was introduced in \cite{BraMos2015} to design polynomial-time algorithms for 
determining $\alpha(G)$. They noted that a largest stable set of a nearly perfect graph can be found in polynomial time. For each vertex $v$, we can find a largest stable set $S_v$ containing $v$ because the non-neihgbours of $v$ induce a perfect graph. Now choose a vertex $a$ that maximizes $S_a$. Then $\{a\} \cup S_a$ is a largest stable set of $G$.   Brandst\"adt and Mosca (\cite{BraMos2015}) proved that (dart, odd hole)-free graphs are nearly perfect. So, there is a polynomial time algorithm
to find a largest stable set of (dart, odd hole)-free graphs.

\section{Conclusions and open problems}\label{sec:conclusions}
In this paper, we study the structure of (dart, odd hole)-free graphs. Our structure results show the existence of a polynomial-time algorithm for  finding a minimum coloring of a (dart, odd hole)-free graph. As expected, our algorithm is reduced to the algorithm to color perfect graphs. The coloring problem for perfect graphs can be expressed as an integer programming  problem. A classic result of Fulkerson (\cite{Ful1973}) provided valuable insights into this problem. His paper 
also suggested the study of the weighted versions of the four optimization problems discussed in Section \ref{sec:optimization}. The algorithm of \cite{GroLov1984} can solve in polynomial time the weighted versions of the four problems.  We define below the weighted version of the coloring problem. 

\noindent {\bf Minimum Weighted Coloring (WCOL)} Given a weighted graph
$G$ such that each vertex $x$ has a weight $w(x)$ which is a
positive integer. Find stable sets $S_1,S_2,\ldots,S_k$ and
integers $I(S_1), \ldots, I(S_k)$ such that for each vertex $x$ we
have $w(x) \leq \Sigma_{x \in S_i}I(S_i)$ and that the sum of the
numbers $I(S_i)$ is minimized. This sum is called the weighted
chromatic number and denoted by $\chi_w(G)$.

Actually, our Theorem~\ref{thm:main} implies the following: if WCOL for graphs with $\alpha = 2$ can be solved in polynomial time, then so can WCOL for   (dart, odd hole)-free graphs. The author (\cite{Hoa2018}) had proposed the following problem.
\begin{problem}
	Determine the complexity of finding a minimum weighted coloring for graphs $G$ with $\alpha(G) = 2$.
\end{problem}
We have shown that (dart, odd hole)-free graphs are perfectly divisible. This result shows $\chi(G) \leq \omega(G)^2$ for a (dart, odd hole)-free graph $G$. We do not know of an odd-hole-free graph $G$ with $\chi(G) = \Omega(\omega(G)^2)$. 
We conclude our paper with the following problem that was first  proposed in \cite{Hoa2018}.
\begin{problem}
	Find an odd-hole-free graph that is not perfectly divisible.
\end{problem}

\end{document}